\documentclass[11pt,a4paper,twoside]{amsart}
\usepackage{amsfonts, amssymb,indentfirst}

\newtheorem{prop}{Proposition}[section]
\newtheorem{theorem}[prop]{Theorem}
\newtheorem{lemma}[prop]{Lemma}

\newtheorem{defi}[prop]{Definition}

\newtheorem{remark}{Remark}
\newenvironment{rem}{\begin{remark}\rm}{\end{remark}}
\newcommand{\cqd}{\hfill$\Box$}

\renewcommand{\hom}[0]{\operatorname{Hom}}
\newcommand{\aut}[0]{\operatorname{Aut}}

\renewcommand{\dim}[0]{\operatorname{dim}}

\newcommand{\G}[0]{$G(2,4)$}

\newcommand{\ch}[0]{\operatorname{ch}}
\newcommand{\Td}[0]{\operatorname{Td}}

\title[Stratification of the Kontsevich moduli space $\overline{M}_{0,n}(\G,d)$]
{ON A STRATIFICATION OF KONTSEVICH'S MODULI SPACE
$\overline{\mathbf{M}}_{\mathbf{0},\mathbf{n}}(\G,\mathbf{d})$ AND
ENUMERATIVE GEOMETRY.}
\author[Cristina Mart{\'\i}nez Ram{\'\i}rez]{ Cristina Mart{\'\i}nez Ramirez}
\address{Max Planck Institute for Mathematics, Vivatsgasse 7, Bonn, 53111.}
 \email{cmartine@mpim-bonn.mpg.de}
\address{Institut for Matematike Fag (CTQM), Aarhus Universitet,
Ny Munkegade, 8000 Aarhus, Denmark} \email{cmartine@imf.au.dk}

\begin{document}
\maketitle

\begin{abstract} We consider a particular
stratification of the moduli space $\overline{M}_{0,n}(G(2,4),d)$ of stable maps to $G(2,4)$.
As an application we compute the
degree of the variety parametrizing rational ruled surfaces with a
minimal directrix of degree $\frac{d}{2}-1$ by studying divisors
in this moduli space of stable maps. For example, there are 128054031872040
rational ruled sextics passing through 25 points in $\mathbb{P}^{3}$ with a
minimal directrix of degree 2.
\end{abstract}

\normalsize \setcounter{equation}{0}

\section{Introduction}

The geometry of  the Kontsevich's moduli space
$\overline{M}_{g,n}(\mathbb{P}^{r},d)$ of degree $d$ stable maps
from $n$-pointed, genus $g$ curves to $\mathbb{P}^{r}$ has been studied
intensively in the literature. R. Vakil studied its connection with the enumerative
geometry of rational and elliptic curves in projective space in
\cite{Vak}. R. Pandharipande studied the theory of
$\mathbb{Q}-$Cartier divisors on this space in \cite{Pan1} and
proved an algorithm for computing all characteristic numbers of
rational curves in $\mathbb{P}^{r}$. R. Pandharipande also
computed the degree of the 1-cuspidal rational locus in the linear
system of degree $d$ plane curves.

The enumerative geometry of rational ruled surfaces of degree $d$
in $\mathbb{P}^{3}$ is closely related to the intrinsic geometry
of the moduli space $\overline{M}_{0,n}(G(2,4),d)$, of Kontsevich stable maps from
$n-$pointed genus 0 curves to the Grassmannian $G(2,4)$ of lines in
$\mathbb{P}^{3}$, representing $d$ times the positive
generator of the homology group $H_{2}(G(2,4),\mathbb{Z})$,
 (see \cite{Mar1}).

In \cite{Mar1} we solved the enumerative problem of
computing the degree of the Severi variety of degree $d$ rational
ruled surfaces in the ambient projective space of surfaces of
degree $d$ in $\mathbb{P}^{3}$. Here we study the more refined
problem of enumerating rational ruled surfaces  with a fixed
directrix of minimum degree. In particular we compute the degree
of the codimension one subvariety of rational ruled surfaces with a minimal directrix of
 degree $\frac{d}{2}-1$ by intersecting divisors in the moduli space of stable maps
$\overline{M}_{0,0}(G(2,4),d)$.

In a first step, we use the variety $R^{0}_{d}$ of degree $d$ morphisms from
$\mathbb{P}^{1}$ to the Grassmannian $G(2,4)$, as a parameter space for
parametrized rational ruled surfaces in $\mathbb{P}^{3}$. We consider a certain
stratification of the variety $R^{0}_{d}$ defined as follows. If
$\mathcal{Q}$ is the universal bundle on $G(2,4)$, the stratum is
defined by the following locally closed condition:

\vspace{0.5cm} $R^{0}_{d,a}=\{f\in R^{0}_{d}| \dim\, H^{1}\,
(\mathbb{P}^{1},f^{*}\mathcal{Q}^{\vee}\otimes
\mathcal{O}_{\mathbb{P}^{1}}(-a-2))\geq 1$ and
$$\dim\,H^{0}(\mathbb{P}^{1},f^{*}_{q}\mathcal{Q}^{\vee}\otimes
\mathcal{O}_{\mathbb{P}^{1}}(-a-1))=0\}.$$

This has a geometrical interpretation as
rational ruled surfaces with a minimal directrix of degree $a$.
First, we consider the closure $R_{d,a}$ of the stratum
$R^{0}_{d,a}$ in the Quot scheme compactification $R_{d}$ of the
space of morphisms $R^{0}_{d}$. We prove that the subschemes
$R_{d,a}$ are irreducible and have the expected codimension as
determinantal varieties (Theorem \ref{determinantal}). We give an
explicit formula for their Poincar\'{e} duals derived from Porteous
formula (Proposition \ref{fundclass}).
In the last section, we consider their closure in the Kontsevich compactification of 
stable maps $\overline{M}_{0,3}(\G,d)$. We study its Picard group (Theorem \ref{generadores}) and
we show that the given stratification of the space $R^{0}_{d}$ can be extended to
$\overline{M}_{0,3}(G(2,4),d)$, but in this case the codimension
of the stratum is different, that is, not all boundary points
correspond to closure points of the open set $R^{0}_{d,a}$, except
for the big stratum of rational ruled surfaces with a minimal
directrix of degree $\frac{d}{2}-1$. In this case we compute its degree by
interpreting it as a tautological intersection on the space of stable maps (Theorem \ref{degree}).
We give explicit enumerative formulas. We show, for example, that there are exactly 128054031872040 rational ruled
sextic surfaces passing through 25 points with a minimal directrix of degree 2.

I. Coskun has approached this problem in
\cite{Cos1} by considering Gromov-Witten invariants for bundles
with a fixed splitting type that correspond to ordinary
Gromov-Witten invariants of certain flag varieties.

 We work over the field of complex numbers $\mathbb{C}$. By a scheme  
 we mean a scheme of finite type over $\mathbb{C}$. By a variety, we
 mean a separated integral scheme. Curves are assumed to be
 complete and reduced. For a variety $X$, $A_{d}X$ and $A^{d}X$
 can be taken to be the Chow homology and
cohomology groups. $A^{d}X$ and $A_{n-d}X$ are identified by the
Poincar\'e duality isomorphism. For $\beta\in A_{k}X$,
$\int_{\beta}c$ is the degree of the zero cycle obtained by
evaluating $c_{k}$ on $\beta$, where $c_{k}$ is the component of
$c$ in $A^{k}X$. If $\mathcal{F}$ is a sheaf over a curve $C$,
$h^{0}(C,\mathcal{F})$ and $h^{1}(C,\mathcal{F})$ denote the
dimensions of the cohomology groups $H^{0}(C,\mathcal{F})$ and
$H^{1}(C,\mathcal{F})$.

\noindent For any formal series
$a(t)=\sum_{k=-\infty}^{k=+\infty}a_{k}\,t^{k}$, we set
$$\triangle_{p,q}(a)=det\,\left( \begin{array}{ccc} a_{p}&\ldots&a_{p+q-1}\\ \vdots&
&\vdots\\ a_{p-q+1}&\ldots & a_{p}
\end{array}\right).$$
\normalsize

\subsection{Acknowledgements}
This paper is part of my Ph.D. thesis. 
I thank the Mathematics Department at Princeton
University, where part of this work was carried out, for their
hospitality and especially professor Rahul Pandharipande,
 for suggestions and comments. I also thank Carolina Araujo,
 Izzet Coskun and Davesh Maulik  for conversations on related
topics. Finally I thank the referee for a long list of comments.

\section{A stratification of the space $R^{0}_{d}$ and enumerative
geometry}
\label{porteous}

Let $R^{0}_{d}$ be the variety of degree $d$ morphisms from
$\mathbb{P}^{1}$ to $G(2,4)$. We will denote $R^{00}_{d}$ the
Zariski open set in $R^{0}_{d}$ corresponding to morphisms which are
birational onto their images.
The Grothendieck Quot scheme $R_{d}$ parametrizing
 rank 2 and degree $d$ quotients of a trivial vector bundle $\mathcal{O}^{4}_{\mathbb{P}^{1}}$,
  is a projective scheme which provides an algebro-geometric compactification of $R^{0}_{d}$. 
  We consider the evaluation map from
$R^{0}_{d}\times \mathbb{P}^{1}$ to the Grassmannian \G, and the
projection map over the first component.

\begin{equation}\label{diagram}\begin{array}{ccccccccc}
 &   & R^{0}_{d}\times \mathbb{P}^{1} & \stackrel{e}{\rightarrow} & G(2,4) & & & & \\
 &   & \pi_{1} \downarrow  & & & & & & \\
  & & R^{0}_{d} & & & & & & \\
\end{array}\end{equation}

The pull back of the universal exact sequence on \G,

\begin{equation}\label{uesg}
0\rightarrow \mathcal{N}\rightarrow
\mathcal{O}_{G(2,4)}^{4}\rightarrow \mathcal{Q}\rightarrow 0
\end{equation}
gives a universal exact sequence on $R^{0}_{d}$ which is the restriction
of the universal quotient over $R_{d}\times \mathbb{P}^{1}$.


 For $f\in R^{0}_{d}$, the pull-back $f^{*}\mathcal{Q}$ under $f$
of the universal quotient on the Grassmannian, is a bundle of rank
2 on $\mathbb{P}^{1}$ isomorphic to
$\mathcal{O}_{\mathbb{P}^{1}}(a)\oplus
\mathcal{O}_{\mathbb{P}^{1}}(d-a)$ where $0\leq a\leq
\frac{d}{2}$, in particular $a\leq d-a$. There is a canonical
exact sequence
$$0\rightarrow N\rightarrow \mathcal{O}^{4}_{\mathbb{P}^{1}}\rightarrow
\mathcal{O}(a)\oplus\mathcal{O}(d-a)\rightarrow 0$$ where $N$ is a
rank 2 and degree $-d$ bundle. Thus, there is a morphism from
$\mathbb{P}(\mathcal{O}(a)\oplus \mathcal{O}(d-a))$ to
$\mathbb{P}^{3}$. The image of this morphism is {\it a scroll}. We
will denote it by $X_{a,d-a}$.

\begin{defi}\label{dim}
For each $0\leq a\leq \frac{d}{2}$, $R^{0}_{d,a}$ is the
subvariety of $R^{0}_{d}$ consisting of those $f\in R^{0}_{d}$
such that $f^{*}\mathcal{Q}\cong
\mathcal{O}_{\mathbb{P}^{1}}(a)\oplus
\mathcal{O}_{\mathbb{P}^{1}}(d-a)$.
\end{defi}

\begin{prop}\label{dim}
 $R^{0}_{d,a}$ is locally closed in $R^{0}_{d}$ and irreducible. If $d>2a$ then it has dimension $3d+5+2a$ and if
 $d=2a$, it has dimension $4d+4$. 
\end{prop}

{\it Proof.} The set $R^{0}_{d,a}$ is locally closed by
Proposition 10 of \cite{Sha}.

\noindent Let
$F_{a}\subset\mathbb{P}\hom(\mathcal{O}^{4},\mathcal{O}_{\mathbb{P}^{1}}(a)\oplus
\mathcal{O}_{\mathbb{P}^{1}}(d-a))$ be the open set of
epimorphisms
$$\mathcal{O}^{4}\twoheadrightarrow
\mathcal{O}_{\mathbb{P}^{1}}(a)\oplus
\mathcal{O}_{\mathbb{P}^{1}}(d-a).$$  $F_{a}$ is in bijection with
the set of locally free rank 2 quotients of
$\mathcal{O}^{4}_{\mathbb{P}^{1}}$ isomorphic to
$\mathcal{O}_{\mathbb{P}^{1}}(a)\oplus
\mathcal{O}_{\mathbb{P}^{1}}(d-a)$ modulo $\mathbb{C}^{*}$. Let
$F_{a}\stackrel{\psi}{\rightarrow} R^{0}_{d,a}$  be the surjective
morphism such that the image by $\psi$ of each quotient is its
class of isomorphism in $R^{0}_{d,a}$. The fibers are isomorphic
to $\mathbb{P}(\aut(\mathcal{O}_{\mathbb{P}^{1}}(a)\oplus
\mathcal{O}_{\mathbb{P}^{1}}(d-a))$.

A dimension calculation yields:

\begin{enumerate}
\item If $a=d-a$,
$$\dim\,\mathbb{P}(\aut(\mathcal{O}(a)\oplus \mathcal{O}(d-a)))=3.$$
\item If $d-a>a$, $$\dim\,\mathbb{P}(\aut(\mathcal{O}(a)\oplus
\mathcal{O}(d-a)))=d-2a+2.$$
\end{enumerate}
The dimension of $F_{a}$ is given by
$$\dim\,\mathbb{P}\hom(\mathcal{O}^{4}_{\mathbb{P}^{1}},
\mathcal{O}_{\mathbb{P}^{1}}(a)\oplus\mathcal{O}_{\mathbb{P}^{1}}(d-a))=\dim\,F_{a}=4d+7.$$
By the fiber dimension theorem (\cite{Sha}, $\S$I.6.3),
\begin{equation}\label{dimd,a}
\dim\,R^{0}_{d,a}=\left\{\begin{array}{ll} 4d+4&\textrm{if $d=2a$}\\
3d+5+2a&\textrm{if $d>2a$}\end{array}\right.
\end{equation}

In particular this proves that the subvarieties $R^{0}_{d,a}$
are irreducible, because $F_{a}$ is irreducible. \cqd

\begin{lemma} \label{stratum}
$q\in R^{0}_{d,a}$ if and only if
$h^{0}(\mathbb{P}^{1},f^{*}_{q}\mathcal{Q}^{\vee}\otimes\mathcal{O}_{\mathbb{P}^{1}}(a))\geq
1$ and
$$h^{0}(\mathbb{P}^{1},f^{*}_{q}\mathcal{Q}^{\vee}\otimes\mathcal{O}_{\mathbb{P}^{1}}(a-1))=0.$$
\end{lemma}

{\it Proof.} Let $q$ be a point in $R^{0}_{d,a}$, and
$f_{q}:\mathbb{P}^{1}\rightarrow G(2,4)$ the morphism
co\-rres\-pon\-ding to that point. Suppose $f_{q}^{*}\mathcal{Q}$
decomposes as $\mathcal{O}_{\mathbb{P}^{1}}(a)\oplus
\mathcal{O}_{\mathbb{P}^{1}}(d-a)$. Tensoring with the line
bundle $\mathcal{O}_{\mathbb{P}^{1}}(a)$,
$$f_{q}^{*}\mathcal{Q}^{\vee}\otimes\mathcal{O}_{\mathbb{P}^{1}}(a)=\mathcal{O}_{\mathbb{P}^{1}}
(2a-d)\oplus\mathcal{O}_{\mathbb{P}^{1}},$$ and tensoring with
$\mathcal{O}_{\mathbb{P}^{1}}(a-1)$, we have
$$f_{q}^{*}\mathcal{Q}^{\vee}\otimes\mathcal{O}_{\mathbb{P}^{1}}(a-1)=\mathcal{O}_{\mathbb{P}^{1}}
(2a-d-1)\oplus\mathcal{O}_{\mathbb{P}^{1}}(-1).$$ Since $a\leq
d-a$, it follows that
$h^{0}(\mathbb{P}^{1},f_{q}^{*}\mathcal{Q}^{\vee}\otimes\mathcal{O}_{\mathbb{P}^{1}}(a))\geq
1,$ and
$$h^{0}(\mathbb{P}^{1},f_{q}^{*}\mathcal{Q}^{\vee}\otimes\mathcal{O}_{\mathbb{P}^{1}}(a-1))=0.$$

Conversely, if $f^{*}_{q}\mathcal{Q}^{\vee}\cong
\mathcal{O}(-n)\oplus \mathcal{O}_{\mathbb{P}^{1}}(n-d)$ with
$n\leq d-n$, then
$$f^{*}_{q}\mathcal{Q}^{\vee}\otimes \mathcal{O}_{\mathbb{P}^{1}}(a)=
\mathcal{O}_{\mathbb{P}^{1}}(n-d+a)\oplus
\mathcal{O}_{\mathbb{P}^{1}}(a-n),$$ 

\noindent then
$h^{0}(\mathbb{P}^{1},f_{q}^{*}\mathcal{Q}^{\vee}\otimes\mathcal{O}_{\mathbb{P}^{1}}(a))\geq
1$ implies $a\geq n$ and

\noindent
$h^{0}(\mathbb{P}^{1},f_{q}^{*}\mathcal{Q}^{\vee}\otimes\mathcal{O}_{\mathbb{P}^{1}}(a-1))=0$
implies $a\leq n$, therefore $n=a$ and $q\in R^{0}_{d,a}$.\cqd

\begin{rem}{\label{dense}}
We observe that the cases $d=2a$ and $d=2a+1$ are special. In
these cases $R^{0}_{d,a}$ has the maximum dimension. This means
that $R^{0}_{d,a}$ is an open set in $R^{0}_{d}$ and therefore the
corresponding degree has already been computed in (\cite{Mar1}).
\end{rem}

\begin{rem}
By the Serre duality theorem and Lemma \ref{stratum} the subvarieties
$R^{0}_{d,a}$ can also be defined as
\begin{equation}\label{defstr}
R^{0}_{d,a}=\{f\in R^{0}_{d}|\, h^{1}\,
(\mathbb{P}^{1},f^{*}\mathcal{Q}\otimes
\mathcal{O}_{\mathbb{P}^{1}}(-a-2))\geq 1 \ \rm{and}\end{equation}
$$ h^{1}\,(\mathbb{P}^{1}, f^{*}\mathcal{Q}\otimes
\mathcal{O}_{\mathbb{P}^{1}}(-(a-1)-2))=0 \}.
$$

\noindent Note that when $d>2a$ a slightly stronger result is true:

\vspace{0.5cm} $q\in R^{0}_{d,a}$ if and only if
$$h^{0}(\mathbb{P}^{1},f^{*}_{q}\mathcal{Q}^{\vee}\otimes\mathcal{O}_{\mathbb{P}^{1}}(a))=1
\ {\rm{and}} \  h^{0}(\mathbb{P}^{1},f^{*}_{q}\mathcal{Q}^{\vee}\otimes\mathcal{O}_{\mathbb{P}^{1}}(a-1))=0.
$$ 

\vspace{0.5cm}
\noindent These sets are locally closed by the Semicontinuity Theorem.

\begin{lemma}The Zariski closure $\overline{R}^{0}_{d,a}$, of the stratum $R^{0}_{d,a}$ in
$R^{0}_{d}$ coincides with
\begin{equation}\label{strata} \{f\in R^{0}_{d}| \, h^{1}\,
(\mathbb{P}^{1},f^{*}\mathcal{Q}\otimes
\mathcal{O}_{\mathbb{P}^{1}}(-a-2))\geq 1 \}.\end{equation}
\end{lemma}
{\it Proof.} By definition, the closure $\overline{R}^{0}_{d,a}$
is the minimum closed set containing $R^{0}_{d,a}$. We observe
that $R^{0}_{d,a}$ is contained in the set given in (\ref{strata}) that
 it is closed, therefore it also contains
 $\overline{R}^{0}_{d,a}$. Now we notice that every closed set
 containing $R^{0}_{d,a}$ contains the closed set

$\{f\in R^{0}_{d}|\, h^{0}(\mathbb{P}^{1},f^{*}_{q}\mathcal{Q}^{\vee}\otimes
\mathcal{O}_{\mathbb{P}^{1}}(a))\geq 1, \ {\rm{and}} \
h^{0}(\mathbb{P}^{1},f^{*}_{q}\mathcal{Q}^{\vee}\otimes
\mathcal{O}_{\mathbb{P}^{1}}(a-1))\geq 0\}$, but the second closed condition
is satisfied trivially, therefore the closure $\overline{R}^{0}_{d,a}$
coincides with the set (\ref{strata}). \cqd

\end{rem}

\begin{rem} We observe that if $f^{*}_{q}\mathcal{Q}\cong \mathcal{O}_{\mathbb{P}^{1}}(-n)\oplus
\mathcal{O}_{\mathbb{P}^{1}}(n-d)$ with $n\leq d-n,$
$h^{0}(\mathbb{P}^{1},f^{*}_{q}\mathcal{Q}^{\vee}\otimes
\mathcal{O}_{\mathbb{P}^{1}}(a))\geq 1$ implies $a\geq n$ and
therefore

\noindent $\overline{R}^{0}_{d,a}=R^{0}_{d,a}\cup R^{0}_{d,a-1}\cup \ldots$
\end{rem}

\begin{lemma}\label{Str1}
The subvariety $R^{0}_{d,a}$ parametrizes rational ruled surfaces
with a minimal directrix of degree $a$.
\end{lemma}
{\it Proof.} The projection
$$\pi_{a}:\mathcal{O}_{\mathbb{P}^{1}}(a)\oplus
\mathcal{O}_{\mathbb{P}^{1}}(d-a)\rightarrow
\mathcal{O}_{\mathbb{P}^{1}}(a),
$$ gives a section $\mathbb{P}(\mathcal{O}_{\mathbb{P}^{1}}(a))
\hookrightarrow\mathbb{P}(f^{*}\mathcal{Q})$. This section is
mapped to a rational curve $C^{a}$ of degree $a$. A study of the
map
\begin{equation}\label{map}
\mathcal{O}_{\mathbb{P}^{1}}(a)\oplus
\mathcal{O}_{\mathbb{P}^{1}}(d-a)\rightarrow
\mathcal{O}_{\mathbb{P}^{1}}(k),
\end{equation}
shows that the map is surjective for $k\geq a$, therefore $C^{a}$
is a {\it directrix} of minimal degree and (\ref{defstr}) has a
geometrical interpretation as the rational ruled surfaces with a
minimal directrix of degree $a$. \cqd

\section{Change of basis in $Pic(R_{d})$.}\label{basis}

In this section we describe geometric
generators for the Picard group of the Quot scheme $R_{d}$.

We consider the cohomology ring of the Grassmannian $G(2,4)$ and
the universal sequence over it:

\begin{equation}\label{uesg}
0\rightarrow \mathcal{N}\rightarrow
\mathcal{O}_{G(2,4)}^{4}\rightarrow \mathcal{Q}\rightarrow 0.
\end{equation}
The special Schubert cycles on $G(2,4)$ can be represented as Chern classes of the
 universal quotient bundle and the universal subbundle:
$$
T_{1}=c_{1}(\mathcal{Q})=c_{1}(\mathcal{N}^{\vee}),
\  T_{b}=c_{2}(\mathcal{N}), \  
\  T_{a}=c_{2}(\mathcal{Q}). \ 
$$

\noindent Also, $T_{3}\in H_{2}(G(2,4),\mathbb{Z})$ will stand for the class
of a line and $T_{4}$ for the class of a point in $G(2,4)$. The
Poincare dual of the hyperplane class $c_{1}(\mathcal{N}^{\vee})$
in $H^{2}(G(2,4),\mathbb{Z})$ determines the Pl\"ucker embedding
of the Grassmannian $G(2,4)$ as a quadric in $\mathbb{P}^{5}$,
that corresponds to a variety of lines in $\mathbb{P}^{3}$. Then
the cycles defined above, have an interpretation as lines in
$\mathbb{P}^{3}$. The cycle $T_{1}$ corresponds to lines in
$\mathbb{P}^{3}$ meeting a given line, $T_{b}$ corresponds to
lines contained in a given plane and $T_{a}$ to lines containing a given point.

\noindent Let $e:R^{0}_{d}\times\mathbb{P}^{1} \rightarrow G(2,4)$
be the evaluation map. The following sets of morphisms define Weil
divisors on $R^{0}_{d}$:
\begin{enumerate}
\item[(A)] The locus of morphisms whose image meets an $a-$plane
$T_{a}$ in the Grassmannian associated to a point $P\in
\mathbb{P}^{3}$. We denote it by:
$$ \label{div1}
D:=\{\varphi \in R^{0}_{d}\ |\ \  e(t,\varphi)\cap T_{a}\neq
\emptyset\}.
$$
\item[(B)] The set of morphisms $Y$ sending a fixed point
$t\in \mathbb{P}^{1}$ to a hyperplane $T_{1}$ on the
Grassmannian:
$$\label{div2} Y:=\{\varphi \in R^{0}_{d}\ |\
e(t,\varphi)\in T_{1}\ \mbox{for  a  fixed } \  t\in
\mathbb{P}^{1} \}.
$$
\end{enumerate}

Since the boundary of $R_{d}$ is of codimension greater than 1
(see Theorem 1.4 of \cite{Ber2}), these divisors extend to divisors in $R_{d}$,
that is $A^{1}(R^{0}_{d})=A^{1}(R_{d})$. Let $Y$, $D$ also denote
the divisors on $R_{d}$. It is clear by the definition that when
we move the point $t\in \mathbb{P}^{1}$ we get a linearly
equivalent divisor. For each index $i\in \mathbb{Z}^{+}$,
we associate a point $t_{i}$ in $\mathbb{P}^{1}$, $Y_{i}$ will denote
the corresponding associated divisor. In the same way, every $a-$plane $T_{a}$
 in $G(2,4)$ is associated to a point $P$ in $\mathbb{P}^{3}$ and if we move the point we also get
 an equivalent divisor with the same index. The divisors $Y_{i}$ were
already considered by Bertram in \cite{Ber1}, (see Lemma 1.1 and Corollary 1.3) where he proved a
moving lemma stating that these varieties can be made to intersect
transversally. S.A. Stromme (Theorem 6.2, \cite{Str}) gives a
basis for the Picard group $A^{1}(R_{d})$, formed by the divisors:

$$\alpha=c_{1}(\pi_{1*}(\mathcal{E}\otimes
\pi_{2}^{*}\mathcal{O}_{\mathbb{P}^{1}}(d))-c_{1}(\pi_{1*}(\mathcal{E}\otimes
\pi_{2}^{*}\mathcal{O}_{\mathbb{P}^{1}}(d-1))),$$
$$\beta=c_{1}(\pi_{1*}(\mathcal{E}\otimes
\pi_{2}^{*}\mathcal{O}_{\mathbb{P}^{1}}(d-1)),$$

\noindent where $\mathcal{E}$ is the universal quotient over
$R_{d}\times \mathbb{P}^{1}$, and $\pi_{1}$, $\pi_{2}$ are the
projections maps over the first and second factors respectively.

\noindent Let $h$ be the positive generator of the Picard group of
$\mathbb{P}^{1}$ as well as its pull-back to $R_{d}\times
\mathbb{P}^{1}$.

\noindent The Chern classes of the universal subbundle
$\mathcal{K}$ are described in the proof of theorem 5.3 of
\cite{Str}, that is
$$A(R_{d}\times \mathbb{P}^{1})=\{A(R_{d})[h] \ |\ h^{2}=0\}.$$ 
Every class $z\in A(R_{d}\times \mathbb{P}^{1})$ can be
written in the form $z=x+h\, y$, with $y=\pi_{1*}(h)\in A(R_{d})$,
and $x=\pi_{1*}(h\,z)\in A(R_{d})$,
\begin{equation}
c_{i}(\mathcal{K})=t_{i}+h\, u_{i-1}\ \ \ (1\leq i\leq 2),
\end{equation}
$t_{i}\in A^{i}(R_{d}),\ u_{i}\in A^{i-1}(R_{d}),\ i=1,2$ as in
\cite{Str}, in particular $u_{0}=-d$.

\vspace{0.5cm} \noindent We express the classes $t_{1},u_{1}$ in
terms of the generators $\alpha$ and $\beta$.

\begin{lemma}The change of basis is
$$\alpha =-t_{1},$$
$$ \beta = u_{1}.$$
\end{lemma}
{\it Proof.} We need to compute the first Chern class of the
bundle,

$$B_{m}=\pi_{1*}(\mathcal{E}\otimes \pi_{2}^{*}\mathcal{O}_{\mathbb{P}^{1}}(m)).$$

\noindent By Riemman Roch, we have that
$$\ch(B_{m})=\pi_{1*}((1+h)\cdot \ch\,(\mathcal{E}\otimes \pi_{2}^{*}\mathcal{O}_{\mathbb{P}^{1}}(m))).$$

\noindent The Chern classes of $\mathcal{E}$ can be computed from
those of $\mathcal{K}$:

$$c_{1}(\mathcal{E})=-c_{1}(\mathcal{K})=-t_{1}+d\,h,$$
$$c_{2}(\mathcal{E})=c^{2}_{1}(\mathcal{K})-c_{2}(\mathcal{K})=
t^{2}_{1}-2d\,h\,t_{1}-t_{2}-u_{1}h.$$

\noindent It follows that,
$$c_{1}(B_{m})=(d-1-m)\,t_{1}+u_{1}.$$

\cqd

\begin{lemma}\label{hypcl} The change of basis from the $\alpha, \beta$ divisors
to the geometric divisors $Y, D$ is:
$$\alpha= Y,$$
$$ \beta=-2dY+D.$$
\end{lemma}

{\it Proof.} By the universal property of the Grassmannian, the
pull back under the evaluation map $e$ of the universal exact
sequence (\ref{uesg}) on $G(2,4)$ gives us a universal exact
sequence on $R^{0}_{d}\times \mathbb{P}^{1}$ which is the
restriction of the universal exact sequence on $R_{d}\times \mathbb{P}^{1}$.
Therefore, by the description of the special Schubert cycles
given, the class [Y] of the divisor $Y$ in the Chow ring $A(R_{d})$
can be identified with the pushforward of the projection map over
the first factor
$\pi_{1*}(he^{*}(T_{1}))=\pi_{1*}(e^{*}(c_{1}(\mathcal{Q}))\cdot
h)=\pi_{1*}(h\,c_{1}(\mathcal{E}))$, where $\mathcal{E}$ is the
universal quotient over $R_{d}\times \mathbb{P}^{1}$. Similarly, the class
of the divisor $D$ can be identified with $\pi_{1*}(c_{2}(\mathcal{E}))$.

The Grothendieck-Riemann-Roch theorem shows that the
geometric divisor $Y$ is just the divisor $\alpha$.  Now, we see
that
$[D]=\pi_{1*}(t_{1}^{2}-2t_{1}\,d\,h-t_{2}+hu_{1})=-2dt_{1}+u_{1}=2d\,\alpha+\beta$,
where $\beta$ is the other generator of the Picard group of
$R_{d}$. In particular, this also proves that the divisors $D$ and $Y$
constitute a basis for the Picard group $A^{1}(R_{d})$. \cqd

\vspace{0.5cm} Let $P_{d}$ be the degree of $R_{d}$ by the
morphism induced by the hyperplane class in $R_{d}$. By lemma
\ref{hypcl}, $P_{d}$ is the degree of the top codimensional
cohomology class given by the self-intersection,
\begin{equation}
P_{d}=\int _{[R_{d}]}[Y]^{4d+4}\cap
[R_{d}]=\int_{[R_{d}]}\alpha^{4d+4}\cap [R_{d}].
\end{equation}

\noindent This intersection number is computed in \cite{RRW} via Quantum
cohomology. It is a Gromov-Witten invariant and can be obtained
too by using the formulas of Vafa and Intriligator,
(see section 5 of \cite{Ber1}).

\noindent In Theorem 4.1 of \cite{Mar1} we prove that the intersection $Y_{1}\cap
Y_{2}\cap Y_{3}\cap D_{1}\cap \ldots \cap D_{4d+1}$ has an excess
intersection component contained in the boundary of $R_{d}$, and therefore the divisors cannot be moved
to make the intersection transversal. We compute it by reinterpreting it as a
Gromov-Witten invariant in the Kontsevich moduli space of stable
maps $\overline{M}_{0,4d+4}(G(2,4),d)$. For $d\geq 3$, this
intersection is $d^{3}\cdot Q_{d}$, where $Q_{d}$ is the Severi
degree of the variety of rational ruled surfaces in
$\mathbb{P}^{3}$, and the factor $d^{3}$ comes from
the multiple covers.

\section{Classes of strata in the Chow ring of $R_{d}$}\label{strata}
In this section we prove that the subschemes $R^{0}_{d,a}$ extend to a
 projective subscheme representing a Chern class and we compute
the classes of the closure of the strata $R^{0}_{d,a}$ in the Quot
scheme compactification $R_{d}$.

\noindent We consider the universal exact sequence over $
R_{d}\times\mathbb{P}^{1}$,
\begin{equation}\label{ues}0\rightarrow \mathcal{K}\rightarrow \mathcal{O}^{4}_{R_{d}\times
\mathbb{P}^{1}}\rightarrow \mathcal{E}\rightarrow 0.\end{equation}

 Tensoring (\ref{ues})  with
$\pi^{*}_{2}\mathcal{O}_{\mathbb{P}^{1}}(-a-2)$, where $\pi_{2}$
is the projection map over the second factor, we obtain the exact
sequence:
$$0\rightarrow \mathcal{K}\otimes
\pi^{*}_{2}\mathcal{O}_{\mathbb{P}^{1}}(-a-2)\rightarrow
\mathcal{O}^{4}_{\mathbb{P}^{1}\times
R_{d}}\otimes\pi^{*}_{2}\mathcal{O}_{\mathbb{P}^{1}}(-a-2)\rightarrow
\mathcal{E}\otimes\pi^{*}_{2}\mathcal{O}_{\mathbb{P}^{1}}(-a-2)\rightarrow
0.$$

 Let $\pi_{1}$ be the projection map over the first factor. The
$\pi_{1*}$ direct image of the above sequence yields an exact
sequence over $R_{d}$:

$$0\rightarrow \pi_{1*}(\mathcal{K}\otimes
\pi^{*}_{2}\mathcal{O}_{\mathbb{P}^{1}}(-a-2))\rightarrow
\pi_{1*}(\mathcal{O}^{4}_{R_{d}\times \mathbb{P}^{1}}\otimes
\pi^{*}_{2}\mathcal{O}_{\mathbb{P}^{1}}(-a-2))\rightarrow $$
$$\rightarrow
\pi_{1*}(\mathcal{E}\otimes\pi^{*}_{2}\mathcal{O}_{\mathbb{P}^{1}}(-a-2))\rightarrow
R^{1}\pi_{1*}(\mathcal{K}\otimes
\pi^{*}_{2}\mathcal{O}_{\mathbb{P}^{1}}(-a-2))\rightarrow $$
$$\rightarrow R^{1}\pi_{1*}(\mathcal{O}^{4}_{\mathbb{P}^{1}\times
R_{d}}\otimes
\pi^{*}_{2}\mathcal{O}_{\mathbb{P}^{1}}(-a-2))\rightarrow
R^{1}\pi_{1*}(\mathcal{E}\otimes\pi^{*}_{2}\mathcal{O}_{\mathbb{P}^{1}}(-a-2))\rightarrow
0. $$

\vspace{0.5cm} By applying the Serre duality theorem, it follows
that
\begin{equation} \label{bundle} R^{1}\pi_{1*}(\mathcal{K}\otimes
\pi^{*}_{2}\mathcal{O}_{\mathbb{P}^{1}}(-a-2))\cong
\pi_{1*}(\mathcal{K}^{\vee}\otimes
\pi_{2}^{*}\mathcal{O}_{\mathbb{P}^{1}}(a))^{\vee}.\end{equation}

We want to see that the sheaf (\ref{bundle}) is a bundle. First,
by the base change theorem, their fibers are isomorphic to
$$H^{1}(\mathbb{P}^{1},\mathcal{K}\otimes \pi_{2}^{*}\mathcal{O}_{\mathbb{P}^{1}}(-a-2)|_{\{q\}\times
\mathbb{P}^{1}})\cong
H^{0}(\mathbb{P}^{1},\mathcal{K}^{\vee}\otimes
\pi_{2}^{*}\mathcal{O}_{\mathbb{P}^{1}}(a)|_{\{q\}\times
\mathbb{P}^{1}})^{\vee}.$$

 It is enough to check 
that the cohomology group
$$H^{1}(\mathbb{P}^{1},\mathcal{K}^{\vee}\otimes
\pi_{2}^{*}\mathcal{O}_{\mathbb{P}^{1}}(a)|_{\{q\}\times
\mathbb{P}^{1}}),$$ vanishes or that the group
$H^{0}(\mathbb{P}^{1},\mathcal{K}\otimes
\pi_{2}^{*}\mathcal{O}_{\mathbb{P}^{1}}(-a-2)|_{\{q\}\times \mathbb{P}^{1}} )$
vanishes, by the Serre duality theorem.

 We can assume that $\mathcal{K}|_{\{q\}\times
\mathbb{P}^{1}}\cong \mathcal{O}_{\mathbb{P}^{1}}(-n)\oplus
\mathcal{O}_{\mathbb{P}^{1}}(n-d)$ with $n\geq 0$ and $n\leq d-n$,
therefore $\mathcal{K}\otimes
\pi_{2}^{*}\mathcal{O}_{\mathbb{P}^{1}}(-a-2)|_{\{q\}\times
\mathbb{P}^{1}}\cong \mathcal{O}_{\mathbb{P}^{1}}(-n-a-2)\oplus
\mathcal{O}_{\mathbb{P}^{1}}(n-d-a-2)$. Since the integers
$-n-a-2$ and $n-d-a-2$ are both negative, the result follows.

For similar reasons, the sheaf
\begin{equation} \label{sheaf2} R^{1}\pi_{1*}(\mathcal{O}^{4}_{R_{d}\times\mathbb{P}^{1}}\otimes
\pi^{*}_{2}\mathcal{O}_{\mathbb{P}^{1}}(-a-2)),\end{equation} is a
bundle over $R_{d}$.

\begin{rem}\label{Eq}
Since the universal quotient sheaf $\mathcal{E}$ is flat over the
Quot scheme $R_{d}$, we have for each $q\in R_{d}$, that
$E_{q}:=\mathcal{E}|_{\{q\}\times \mathbb{P}^{1}}$ is a coherent
sheaf over $\mathbb{P}^{1}$ and therefore it can be written as a direct
sum of its locally free part plus its torsion. By the Riemann-Roch
theorem, we have that
\begin{equation}\label{RR}
h^{0}(E_{q})-h^{1}(E_{q})=deg\,(E_{q})+rank \,(E_{q})(1-g)=d+2.
\end{equation}
In particular, the degree of $E_{q}$ is constant for every $q\in R_{d}$. The sets
$$h^{1}(E_{q})\geq r,$$ where $r$ is an integer, are closed by the Semicontinuity
theorem and by the identity (\ref{RR}), they can also be defined as:
$$h^{0}(E_{q})\geq d+2+r.$$
\end{rem}

Let us consider the Zariski closure $R_{d,a}$ of the sets
$R^{0}_{d,a}$ inside the Quot scheme compactification of the space of morphisms:
$$\{q\in R_{d}|\, h^{0}(\mathbb{P}^{1},E^{\vee}_{q}\otimes \mathcal{O}_{\mathbb{P}^{1}}(a))\geq
1\},$$ or equivalently by Serre duality:

$$\{q\in R_{d}|\, h^{1}(\mathbb{P}^{1},E_{q}\otimes \mathcal{O}_{\mathbb{P}^{1}}(-a-2))\geq
1\}.$$

\begin{theorem} \label{determinantal}
$R_{d,a}$ is the locus where the map
\begin{equation}\label{conucleo}R^{1}\pi_{1*}(\mathcal{K}\otimes
\pi^{*}_{2}\mathcal{O}_{\mathbb{P}^{1}}(-a-2))\rightarrow
R^{1}\pi_{1*}(\mathcal{O}^{4}_{R_{d}\times\mathbb{P}^{1} }\otimes
\pi^{*}_{2}\mathcal{O}_{\mathbb{P}^{1}}(-a-2)).\end{equation} is
not surjective. It is non-empty and has the expected codimension $d-2a-1$
as a determinantal variety.
\end{theorem}

{\it Proof.} The map (\ref{conucleo}) is not surjective in the
support of the sheaf
\begin{equation}\label{sheaf}
R^{1}\pi_{1*}(R_{d}\times\mathbb{P}^{1}
,\mathcal{E}\otimes\pi^{*}_{2}\mathcal{O}_{\mathbb{P}^{1}}(-a-2)),
\end{equation}
that is, in the points $q\in R_{d}$ such that
$$h^{1}\,(\mathbb{P}^{1}, \ E_{q}\otimes
\mathcal{O}_{\mathbb{P}^{1}}(-a-2))\geq 1,$$ or equivalently by Serre duality in
$R_{d,a}$, or by the observation in Remark \ref{Eq} in
$$\{q\in R_{C,d}\,|h^{0}\,(\mathbb{P}^{1}, \ E_{q}\otimes \mathcal{O}_{\mathbb{P}^{1}}(-a-2))\geq d-2a-1\}.$$

The non-emptiness follows from the fact that given a vector
 bundle $E\cong \mathcal{O}_{\mathbb{P}^{1}}(a)\oplus
 \mathcal{O}_{\mathbb{P}^{1}}(d-a)$, it is generated by global sections since by assumption $a,\ d-a\geq
 0$ and therefore there exists
 $f\in R^{0}_{d}$ such that $f^{*}\mathcal{Q}\cong E$.

By applying the Serre duality theorem we get that

\begin{equation}
R^{1}\pi_{1*}(\mathcal{K}\otimes\pi^{*}_{2}\mathcal{O}_{\mathbb{P}^{1}}(-a-2))\cong
\pi_{1*}(\mathcal{K}^{\vee}\otimes\pi_{2}^{*}\mathcal{O}_{\mathbb{P}^{1}}(a))^{\vee}.
\end{equation}


We apply Riemann-Roch in $\mathbb{P}^{1}$ in order to compute the
dimension of their fibers:
$$\dim\,H^{0}(\mathbb{P}^{1} , \mathcal{K}^{\vee}\otimes
\pi_{2}^{*}\mathcal{O}_{\mathbb{P}^{1}}(a)|_{\{p\}\times
\mathbb{P}^{1}})=d+2a+2,\ p\in R_{d}.$$ As a consequence
$\pi_{1*}(\mathcal{K}^{\vee}\otimes
\pi_{2}^{*}\mathcal{O}_{\mathbb{P}^{1}}(a))$ is a bundle of rank
$d+2a+2$ over $R_{d}$. Again we see by the Serre duality theorem that

 \begin{equation} R^{1}\pi_{1*}(\mathcal{O}^{4}_{R_{d}\times
\mathbb{P}^{1}}\otimes\pi_{2}^{*}\mathcal{O}(-a-2)) \cong
\pi_{1*}(\mathcal{O}^{4}_{R_{d}\times \mathbb{P}^{1}}\otimes
\pi_{2}^{*}\mathcal{O}(-a-2))^{\vee}.
\end{equation}



Since $\mathcal{O}^{4}_{\mathbb{P}^{1}}\otimes
\pi_{2}^{*}\mathcal{O}_{\mathbb{P}^{1}}(a)|_{\{p\}\times
\mathbb{P}^{1}}$ is a locally free sheaf over $\mathbb{P}^{1}$, it
can be decomposed over $\mathbb{P}^{1}$, so that
$$\dim\, H^{0}(\mathbb{P}^{1}
,\mathcal{O}^{4\vee}_{R_{d}\times
\mathbb{P}^{1}}\otimes\pi_{2}^{*}\mathcal{O}_{\mathbb{P}^{1}}(a)|_{\mathbb{P}^{1}\times\{p\}})=4(a+1).$$
As a consequence $\pi_{1*}(\mathcal{O}^{4}_{R_{d}\times
\mathbb{P}^{1}}\otimes
\pi_{2}^{*}\mathcal{O}_{\mathbb{P}^{1}}(a))$ is a trivial bundle
of rank $4a+4$ over $R_{d}$, that is,
$$\pi_{1*}(\mathcal{O}^{4}_{R_{d}\times \mathbb{P}^{1}}\otimes
\pi_{2}^{*}\mathcal{O}_{\mathbb{P}^{1}}(a))\cong
\mathcal{O}_{R_{d}}^{4a+4}.$$

Finally, all these facts give us the following morphism of
bundles:
$$\pi_{1*}(\mathcal{K}^{\vee}\otimes
\pi_{2}^{*}\mathcal{O}_{\mathbb{P}^{1}}(a))\stackrel{\phi}{\rightarrow}\mathcal{O}^{4a+4}_{R_{d}}.$$

The expected codimension of $R_{d,a}$ as a determinantal variety
is: $$((d+2a+2)-(4a+3))\cdot((4a+4)-(4a+3))=d-2a-1\ (\textrm{$\S$
II.4, \cite{ACGH}}).$$ This codimension coincides with the
computed codimension of $R^{0}_{d,a}$ in \ref{dim}. 
\cqd

\subsubsection{Class of $R_{d,a}$ in $A(R_{d})$.}
\begin{lemma}\label{irreducible}
$R_{d,a}$ is irreducible and does not have components contained in the
boundary $R_{d,a}-R^{0}_{d,a}$.
\end{lemma}
{\it Proof.} Since the subschemes $R_{d,a}$ are the closure of the
irreducible sets $R^{0}_{d,a}$, it follows that they are also
irreducible, and consequently cannot have components contained at infinity. \cqd

\begin{prop}\label{fundclass} If $R_{d,a}$ is either empty or has the expected
codimension $d-2a-1$, its fundamental class is given by the formula:
   $$[R_{d,a}]=-c_{d-2a-1}(\pi_{1*}(\mathcal{K}^{\vee}
\otimes \pi_{2}^{*}(\mathcal{O}_{\mathbb{P}^{1}}(a)))).$$
\end{prop}

{\it Proof.} By Lemma \ref{irreducible} and Theorem \ref{determinantal}, the subschemes $R_{d,a}$ are irreducible and have the
expected codimension as determinantal varieties . The formula for their Poincare duals
is derived from Porteous formula:
$$ [R_{d,a}]=\triangle_{d-2a-1,1}(c_{t}(-\pi_{1*}(\mathcal{K}^{\vee}\otimes
\pi_{2}^{*}(\mathcal{O}_{\mathbb{P}^{1}}(a))))) $$

\cqd

\subsubsection*{Computation of Chern classes.}
In order to compute the Chern classes of the bundle
$\pi_{1*}(\mathcal{K}^{\vee}\otimes
\pi_{2}^{*}\mathcal{O}_{\mathbb{P}^{1}}(a))$ we apply
the Grothendieck-Riemann-Roch theorem, so that
\begin{equation}\label{RRG}
\ch\,(\pi_{1*}(\mathcal{K}^{\vee}\otimes
\pi_{2}^{*}(\mathcal{O}_{\mathbb{P}^{1}}(a))))=\pi_{1*}(\Td(R_{d}\times
\mathbb{P}^{1}/R_{d})\cdot \ch(\mathcal{K}^{\vee}\otimes
\pi_{2}^{*}\mathcal{O}_{\mathbb{P}^{1}}(a))).
\end{equation}

\begin{equation}
c_{i}(\mathcal{K})=t_{i}+h\, u_{i-1}\ \ \ (1\leq i\leq 2),
\end{equation}

\noindent $t_{i}, u_{i}\in A^{i}(R_{d}),i=1,\ldots,n-1$ as in
\cite{Str}, in particular $u_{0}=-d$.

\vspace{0.5cm} \noindent First we compute the Chern classes of
$\mathcal{K}^{\vee}\otimes
\pi_{2}^{*}\mathcal{O}_{\mathbb{P}^{1}}(a)$:
\begin{align*}
c_{1}(\mathcal{K}^{\vee}\otimes
\pi_{2}^{*}\mathcal{O}_{\mathbb{P}^{1}}(a))&=c_{1}(\mathcal{K}^{\vee})+
c_{1}(\pi_{2}^{*} \mathcal{O}_{\mathbb{P}^{1}}(a))\cdot
c_{0}(\mathcal{K}^{\vee})\\
&=h\,(2a+d)-t_{1},
\end{align*}
\begin{align*}c_{2}(\mathcal{K}^{\vee}\otimes \pi_{2}^{*}\mathcal{O}_{\mathbb{P}^{1}}(a))&=
c_{2}(\mathcal{K}^{\vee})+ c_{1}(\mathcal{K}^{\vee}) \cdot
c_{1}(\pi_{2}^{*}\mathcal{O}_{\mathbb{P}^{1}}(a))+c_{1}(\pi_{2}^{*}\mathcal{O}_{\mathbb{P}^{1}}(a))^2
\cdot c_{0}(\mathcal{K}^{\vee})\\ &= -a\, t_{1}\, h-t_{2}-h\,
u_{1},
\end{align*}
\begin{align*}c_{3}(\mathcal{K}^{\vee}\otimes \pi_{2}^{*}\mathcal{O}_{\mathbb{P}^{1}}(a))=
a\,h t_{2}.
\end{align*}
 Let us denote by ${\rm{ch}}_{i}$, the $i-$homogeneous
part of the Chern character of a bundle.

$$\ch_{0}\,(\mathcal{K}^{\vee}\otimes
\pi_{2}^{*}\mathcal{O}_{\mathbb{P}^{1}}(a)))=2,$$

$$\ch_{1}\,(\mathcal{K}^{\vee}\otimes
\pi_{2}^{*}\mathcal{O}_{\mathbb{P}^{1}}(a)))=-t_{1}+h\,(d+2a), $$

\begin{align*}\ch_{2}\,(\mathcal{K}^{\vee}\otimes \pi_{2}^{*}\mathcal{O}_{\mathbb{P}^{1}}(a))) &=\frac{1}{2}
[t_{1}^{2}-2t_{1}\,h\,(2a+d)+2at_{1}h+2t_{2}+2hu_{1}
\end{align*}

\begin{align*}  -3u_{2}\,h-3a\,h\,t_{2}],
\end{align*}

\begin{lemma}
$$\ch_{n}\,(\mathcal{K}^{\vee}\otimes \pi_{2}^{*}\mathcal{O}_{\mathbb{P}^{1}}(a))=
{\rm{coeff}}_{(t^{n})}\left(\sum_{n}\frac{(-1)^{n-1}}{n}
\sum_{j}c_{j}(\mathcal{K}^{\vee}\otimes
\pi_{2}^{*}\mathcal{O}_{\mathbb{P}^{1}}(a)))t^{j}\right)^{n}.
$$
\end{lemma}

{\it Proof.} The Chern character of a bundle $E$ is a polynomial
in the Chern classes $x_{i}=c_{i}(E)$ defined by the formula:
$$\ch(E)=\sum_{i}e^{x_{i}}=\sum_{i}\sum_{n}\frac{(x_{i})^{n}}{n!}\in H^{*}(X,\mathbb{Q}).$$
If we call $\sigma_{i}$ to the symmetric functions, that is,

$$\sum_{r=0}^{n}\sigma_{r}t^{r}=\prod_{i=1}^{n}(1+x_{i}t),$$ we
get that $$log\,
(1+\sigma_{1}t+\sigma_{2}t^{2}+\ldots)=log\prod_{i}(1+x_{i}t)$$

$$=\sum_{i}\sum_{n}
\frac{(-1)^{n-1}}{n}(x_{i}t)^{n}=\sum_{n}\frac{(-1)^{n-1}}{n}(\sum
x_{i}^{n})t^{n}.$$ and from here the formula above follows.
\cqd

\noindent Now applying ($\ref{RRG}$), we have
$$\ch(\pi_{1*}(\mathcal{K}^{\vee}\otimes \pi^{*}\mathcal{O}_{\mathbb{P}^{1}}(a)))=
\pi_{1*}\left( (1+h)\,(\ch\,(\mathcal{K}^{\vee}\otimes
\mathcal{O}_{\mathbb{P}^{1}}(a)))\right ).$$

Therefore the Chern classes of $\pi_{1*}(\mathcal{K}^{\vee}\otimes
\pi_{2}^{*}\mathcal{O}_{\mathbb{P}^{1}}(a))$ are:
\begin{align*}c_{1}(\pi_{1*}(\mathcal{K}^{\vee}\otimes
\pi_{2}^{*}\mathcal{O}_{\mathbb{P}^{1}}(a)))=-t_{1}\,(a+d+1)+u_{1},\end{align*}
$$c_{2}(\pi_{1*}(\mathcal{K}^{\vee}\otimes
\pi_{2}^{*}\mathcal{O}_{\mathbb{P}^{1}}(a)))=\frac{1}{2}\,c_{1}^{2}-\left[
\frac{1}{2}\,(2a+d)\,t_{1}^{2}-\frac{1}{2}\,a\,t_{1}^{2}+\frac{1}{2}\,
(2a+d)\,t_{2}-\frac{1}{2}u_{1}t_{1}\right],$$
\begin{center}
$\vdots$
\end{center}
\begin{lemma}
$$c_{n}(\pi_{1*}(\mathcal{K}^{\vee}\otimes
\pi_{2}^{*}\mathcal{O}_{\mathbb{P}^{1}}(a)))=$$
$$-\sum_{r=1}^{n}\frac{(-1)^{r-1}}{n}r!\ch_{r}(\pi_{1*}(\mathcal{K}^{\vee}\otimes
\pi_{2}^{*}\mathcal{O}_{\mathbb{P}^{1}}(a)))\,c_{n-r}(\pi_{1*}(\mathcal{K}^{\vee}\otimes
\pi_{2}^{*}\mathcal{O}_{\mathbb{P}^{1}}(a))).$$
\end{lemma}

{\it Proof.} For each $r\geq 1$ the $r'th$ power sum is:
$$p_{r}=\sum x_{i}^{r}=m_{(r)}.$$ The generating function for the
$p_{r}$ is: \begin{equation}\label{eq1} p(t)=\sum_{r\geq
1}p_{r}t^{r-1}=\sum_{i\geq 1}\sum_{r\geq 1}x_{i}^{r}t^{r-1}=
\sum_{i\geq 1}\frac{d}{dt}log\frac{1}{1-x_{i}t}.\end{equation}
\begin{equation}\label{eq2}P(t)=\frac{d}{dt}\prod_{i\geq
1}(1-x_{i}t)^{-1}=\frac{d}{dt}log H(t)=\frac{H'(t)}{H(t)}$$
$$P(-t)=\frac{d}{dt}log\, E(t)=\frac{E'(t)}{E(t)}.\end{equation}
From (\ref{eq1}) and (\ref{eq2}) we get that
$$n\sigma_{n}=\sum_{r=1}^{n}(-1)^{r-1}p_{r}\sigma_{n-r},$$ and the
formula above follows. \cqd

\section{Another compactification of the space $R^{0}_{d}$.}

Now we will consider the coarse moduli space $\overline{M}_{0,n}(G(2,4),d)$
of degree $d$ Kontsevich stable maps from $n$-pointed, genus 0 curves to
$G(2,4)$.
It is smooth, has dimension $4d+1+n$ and there is a universal
family over it:

\vspace{0.5cm}
\begin{picture}(200,80)
\put(165,70){$0 \rightarrow \mathcal{N} \rightarrow
\mathcal{O}^{4}\rightarrow  \mathcal{Q} \rightarrow  0$}
\put(222,55){$\downarrow$}

\put(50,40){$\overline{M}_{0,n+1}(G(2,4), d)$}
\put(136,42){\vector(1,0){60}} \put(155,45){$e_{i}$}
\put(200,42){$G(2,4)$} \put(65,28){$\pi_{i} \downarrow  \uparrow
\sigma_{i}$} \put(50,10){ $\overline{M}_{0,n}(G(2,4), d)$ }
\end{picture}

\noindent Here $\pi_{n}:\overline{M}_{0,n+1}(G(2,4), d)\rightarrow
\overline{M}_{0,n}(G(2,4), d)$ is the forgetful morphism, which
consists of forgetting the mark $p_{1}$ and $\sigma_{i}$ is the
section corresponding to the mark $p_{i}$. The image of
$\sigma_{i}$ is the (closure of the) locus of maps whose source
curve has two components, one of which is rational, carries just
the two marks $p_{i}$ and $p_{n+1}$ and is contracted by $\mu$.
For $n=3$, the space $\overline{M}:=\overline{M}_{0,3}(G(2,4),d)$ provides a
compactification of the space $R^{0}_{d}$.
The marked points yield
canonical line bundles $\mathcal{L}_{i}=e_{i}^{*}(T_{1})$ on
$\overline{M}$ via the $i$-evaluation map, where $T_{1}$ is the
hyperplane class in $G(2,4)$. Since the boundary is of pure
codimension 1 in $\overline{M}_{0,n}(G(2,4),d)$, each irreducible component is a
Weil divisor. Let $\triangle$ be the set of components of the
boundary. For $d\geq 1$, Weil divisors are obtained on
$\overline{M}_{0,n}(G(2,4),d)$ by considering the locus 
corresponding to maps meeting an $a-$plane and a $b-$plane on the
Grassmannian. We denote the corresponding classes in
$Pic(\overline{M})\otimes \mathbb{Q}$ by $\mathcal{A}$ and
$\mathcal{B}$ respectively. Let $w_{\pi}$ denote the relative
dualizing sheaf, then a map $\mu \in \overline{M}_{0,n}(G(2,4),
d)$ is stable if and only if $w_{\pi}(p_{1}+\ldots +p_{n})\otimes
\mu^{*}(T_{3})$ is ample. In particular
$deg\,(w_{\pi}(p_{1}+\ldots +p_{n})\otimes \mu^{*}(T_{3}))>0$, so
that if $E\subset C$ such that $\dim \mu(E)=0$, then
$deg\,(w_{\pi}(p_{1}+\ldots +p_{n}))>0$.

\begin{theorem}\label{generadores}
If $d\geq 1$, then $Pic(\overline{M})$ is generated by
$\triangle\cup \{\mathcal{A}\}\cup \{\mathcal{L}_{i}\}_{i=1,2,3}$.
\end{theorem}
{\it Proof.} Consider the universal family of stable degree $d$
maps of 3-pointed curves $R^{00}_{d}\times
\mathbb{P}^{1}\rightarrow G(2,4)$. By the universal property,
there is an injection
$$R^{00}_{d}\hookrightarrow \overline{M}.$$ The complement of
$R^{00}_{d}$ is the boundary of $\overline{N}$ and it is also the locus
co\-rres\-pon\-ding to maps with reducible domain. Since the
complement of $R^{00}_{d}$ in $R^{0}_{d}$ is of codimension
greater or equal to 2 (see Lemma 2.5 of \cite{Mar1}), the Picard
group of $R^{00}_{d}$ coincides with the Picard group of
$R^{0}_{d}$. The divisor $\mathcal{A}$ restricted to $R^{0}_{d}$,
i.e. to the moduli points corresponding to stable maps with
irreducible source, is linearly equivalent to $D_{i}$ and
$\mathcal{L}_{i}$ is li\-near\-ly equivalently to the
corresponding divisor $Y_{i}$. In these points, markings are not
necessary since the maps are already stable. In section
\ref{basis}, we argued that the divisors $Y_{i}$ and $D_{i}$
constitute a basis of $A^{1}(R^{0}_{d})$, so that
$A_{4d+3}(\overline{M})$ is generated by the boundary $\triangle$
and the Weil divisors $\mathcal{L}_{i}$ ($i=1,2,3$), and
$\mathcal{A}$. We note that since $Y_{1}\sim Y_{2}\sim Y_{3}$ in
$R^{0}_{d}$, each one of the $\mathcal{L}_{i}$ generates the
Picard group.\cqd

\vspace{0.5cm} \noindent If $d=0$, then
$\overline{M}_{0,n}(G(2,4),0)\cong \overline{M}_{0,n}\times
G(2,4)$ and $\mathcal{L}_{i}$ is the pull-back of
$\mathcal{O}_{G(2,4)}(1)$ from the second factor.

\subsubsection*{Boundary of $\overline{M}_{0,n}(G(2,4),d)$}
The irreducible components of the boundary are in bijective
correspondence with the data of weighted partitions $(A\cup B,
d_{A},d_{B})$, where
\begin{itemize}
\item $A\cup B$ is a partition of the set of marked points. There
are only 2 possible partitions: 1+2, 3+0. \item $d_{A}+d_{B}=d$,
$d_{A}\geq 0, d_{B}\geq 0$, \item if $d_{A}=0$ (resp. $d_{B})$,
then $|A|\geq 2$ (resp. $|B|\geq 2$).
\end{itemize}

\noindent A component $(A\cup B, d_{A},d_{B})$ consists
generically of maps with a union of $\mathbb{P}^{1}$'s as domains.
Each one gives a curve in $G(2,4)$ of degrees $d_{A}$ and $d_{B}$,
respectively.

\subsection{Stratification of $\overline{M}$.}

The stratification of the space $R^{0}_{d}$ given in \ref{Str1}
can be extended to $\overline{M}$, we can
even compute the classes of the strata in the Chow ring
$A(\overline{M})$, but they have no enumerative meaning.

\subsubsection{Closure of the stratum in $\overline{M}$.}\label{coskun} The
closure of the stratum has already been studied by I. Coskun in
\cite{Cos1} by taking flat limits of surfaces in one parameter
families.

Let $V_{d,a}$ denote the closure of $R^{0}_{d,a}$ in
$\overline{M}$. Let $K^{i}$ be the boundary component
corresponding to the degree partition $i+(d-i)=d$. Geometrically
these boundary points correspond to reducible curves in the
Grassmannian $G(2,4)$ of degree $i$, $d-i$ respectively, or
equivalently a union of irreducible rational ruled surfaces of
degree $i$ and $d-i$ in $\mathbb{P}^{3}$.

\noindent A generic intersection point of the boundary with $V_{d,a}$ corresponds
to a moduli point, $[C, p_{1},p_{2},p_{3},\mu:C\rightarrow G(2,4)]$, where
$C=C_{1}\cup C_{2}$ is a union of two $\mathbb{P}^{1}$'s with
$\mu^{*}\mathcal{Q}|_{C_{1}}\cong
\mathcal{O}_{\mathbb{P}^{1}}(a_{1})\oplus
\mathcal{O}_{\mathbb{P}^{1}}(b_{1})$, $a_{1}+b_{1}=i$ and
$\mu^{*}\mathcal{Q}|_{C_{2}}\cong
\mathcal{O}_{\mathbb{P}^{1}}(a_{2})\oplus
\mathcal{O}_{\mathbb{P}^{1}}(b_{2})$, $a_{2}+b_{2}=d-i$. The
integers $a_{j},b_{j}$ satisfy $a_{j}\geq 0$, $b_{j}\geq 0 $,
$\sum_{j}(a_{j}+b_{j})=d$ and $a_{1}+a_{2}\leq a$. This last
condition means that the flat limit of the directrices is a
connected curve of degree $a$ whose restriction to each of the
surfaces is in a section class. Conversely, given any connected
curve $C$ of degree $a\leq \sum\,(a_{j}+b_{j})/2$ whose
restriction to each component is in a section class, there exists
a one parameter family of $X_{a,b}$ specializing to the reducible
surface such that the limit of the directrices is $C$.

\vspace{0.5cm} The stratification given in section \ref{porteous}
is extended to $\overline{M}$ by defining the stratum as
\begin{equation}\label{str2}\overline{M}_{d,a}:= \{\mu\in \overline{M}|
h^{1} ( \mu^{*}\mathcal{Q}\otimes w_{\pi}(-a\,p_{i}))\geq 1\}.
\end{equation}
This set is a closed subset of $\overline{M}$ by the
Semicontinuity Theorem.

\begin{prop}
$\overline{M}_{d,a}$ is the locus where the map
\begin{equation}\label{coker2}R^{1}\pi_{i*}(e_{i}^{*}(\mathcal{N})\otimes
w_{\pi_{i*}}(-a\,p_{i}))\rightarrow
R^{1}\pi_{i*}(e_{i}^{*}(\mathcal{O}^{4})\otimes
w_{\pi_{i*}}(-a\,p_{i}))\end{equation} is not surjective.
\end{prop}
{\it Proof.} We proceed as in section \ref{strata}, considering
the long exact sequence:

$$0\rightarrow \pi_{i*}(e_{i}^{*}\mathcal{N}\otimes w_{\pi}(-a\,p_{i}))\rightarrow
\pi_{i*}(e_{i}^{*}\mathcal{O}^{4}\otimes
w_{\pi}(-ap_{i}))\rightarrow \pi_{*}(e_{i}^{*}\mathcal{Q}\otimes
w_{\pi}(-ap_{i}))$$ $$\rightarrow
R^{1}\pi_{*}(e_{i}^{*}(\mathcal{N})\otimes
w_{\pi}(-a\,p_{i}))\rightarrow
R^{1}\pi_{*}(e_{i}^{*}(\mathcal{O}^{4})\otimes
w_{\pi}(-a\,p_{i}))\rightarrow
$$ $$ R^{1}\pi_{*}(e_{i}^{*}(\mathcal{Q})\otimes w_{\pi}(-a\,p_{i}))\rightarrow 0.$$

\noindent The map (\ref{coker2}) is not surjective in the
support of the sheaf $ R^{1}\pi_{1*}(e_{i}^{*}(\mathcal{Q})\otimes
w_{\pi}(-a\,p_{i}))$, that is at the moduli points
$[C,p_{1},p_{2},p_{3},\mu:C\rightarrow G(2,4)]$ where $h^{1}(C,
\mu^{*}\mathcal{Q}\otimes w_{\pi}(-a\,p_{i}))\geq 1,$ and these
are the points defining the stratum (\ref{str2}).\cqd

\subsubsection{Enumerative geometry of $\overline{M}_{d,a}$.}
\vspace{0.5cm} We would like to know which points in the boundary
of $\overline{M}$ satisfy the condition:
$$\label{condition}h^{1}(C, \mu^{*}\mathcal{Q}\otimes
w_{\pi}(-a\,p_{i}))\geq 1.$$ The next theorem will address this
question.

Recalling the notation of \cite{Mar1}, $Q_{d}$ will be the number
of rational ruled surfaces through $4d+1$ points in
$\mathbb{P}^{3}$. 
The number $Q_{d}$ is the zero-cycle $\mathcal{A}^{4d+1}$ in the Chow ring of
$\overline{M}_{0,0}(G(2,4),d)$. The number $Q_{d}^{b}$ will denote
the zero-cycle $\mathcal{A}^{4d}\mathcal{B}$, or equivalently the
Gromov-Witten invariant,
$I_{0,4d+1,d}(T_{a},\stackrel{4d}{\ldots},T_{a},T_{b})$. The
following table lists the numbers $Q_{d}^{b}$ for each degree
$1\leq d\leq 9$. To compute these numbers we have applied Farsta,
a program due to Andrew Kresch. These numbers correspond
geometrically to rational ruled surfaces of degree $d$ through
$4d$ points in $\mathbb{P}^{3}$ and tangent to a fixed plane.

\vspace{0.5cm}
\begin{center}
\begin{tabular}{|c|r|r|}
\hline $d$ & $Q^{b}_{d}$  & $4d$   \\
\hline
 1 & 0& 4 \\
 2 & 2& 8 \\
 3 & 1824 & 12 \\
 4 & 3094440 & 16 \\
 5 & 15383867920 & 20 \\
 6 & 188115939619440 & 24 \\

 7 & 3000224401806629008 & 28 \\
 8 & 219761533783440334862592 & 32 \\
 9 & 17394248462381072210049044320 & 36 \\
 \hline
 \end{tabular}
 \end{center}

\begin{theorem}\label{degree} $V_{d,a}\subset \overline{M}_{d,a}$ and $V_{d,a}=\overline{M}_{d,a}$
only when $a=\frac{d}{2}-1$. In this case, the stratum
$R^{0}_{d,a}$ is of pure codimension 1 and its degree is given by
the formula
$$\frac{1}{2}\,d^{3}\cdot (Q_{d}+ Q_{d}^{b}).$$
\end{theorem}

{\it Proof.} By semicontinuity (see \cite{Har}, III. 12.8), the
points

\noindent $[C,p_{1},p_{2},p_{3},\mu:C\rightarrow G(2,4)]$ in the closure of
$R^{0}_{d,a}$ in $\overline{M}$ satisfy the condition
\ref{condition}. The converse is false. We will construct some families satisfying the
condition but which are not in the closure of the stratum.

Grothendieck duality implies $$h^{1}(C,\mu^{*}\mathcal{Q}\otimes w_{\pi}(-a\,p_{i}))=h^{0}(C,
\mu^{*}\mathcal{Q}^{\vee}\otimes \mathcal{O}(ap_{i})).$$ For a
general boundary point, the source curve $C=C_{1}\cup C_{2}$ is a
nodal curve. There is a canonical exact sequence:
$$0\rightarrow \mu^{*}\mathcal{Q}^{\vee}|_{C}\otimes
\mathcal{O}(ap_{i}) \rightarrow
\mu^{*}\mathcal{Q}^{\vee}|_{C_{1}}\otimes \mathcal{O}(ap_{i})
\oplus \mu^{*}\mathcal{Q}^{\vee}|_{C_{2}}\otimes
\mathcal{O}(ap_{i})\rightarrow \mathbb{C}^{2}\rightarrow 0.$$ This
gives the following injections:

\begin{picture}(200,60)
\put(0,40){$H^{0}(C_{1}, \mu^{*}\mathcal{Q}^{\vee}|_{C_{1}}\otimes
\mathcal{O}(ap_{i}))\rightarrow
H^{0}(C,\mu^{*}\mathcal{Q}^{\vee}|_{C}\otimes
\mathcal{O}(ap_{i}))$} \put(190,20){$\uparrow$}
\put(160,0){$H^{0}(C_{2},\mu^{*}\mathcal{Q}^{\vee}|_{C_{2}}\otimes
\mathcal{O}(ap_{i}))$.}
\end{picture}

\vspace{0.5cm} \noindent Since $C_{i}\cong \mathbb{P}^{1}$, it can
be assumed that $\mu^{*}\mathcal{Q}|_{C_{1}}\cong
\mathcal{O}(a_{1})\oplus \mathcal{O}(b_{1})$ and
$\mu^{*}\mathcal{Q}|_{C_{2}}\cong \mathcal{O}(a_{2})\oplus
\mathcal{O}(b_{2})$ with $a_{1},a_{2},b_{1},b_{2}$ being integers
greater than or equal to 0 and satisfying $\sum_{i}a_{i}+b_{i}\leq d$.
Also we assume that $P_{i}\in C_{1}$, so that $h^{0}(C_{2},
\mu^{*}\mathcal{Q}^{\vee}|_{C_{2}}\otimes
\mathcal{O}(ap_{i}))=h^{0}(C_{2},\mathcal{O}(-a_{2})\oplus
\mathcal{O}(-b_{2}))$, and

\noindent $h^{0}(C_{1}, \mu^{*}\mathcal{Q}^{\vee}|_{C_{1}}\otimes
\mathcal{O}(ap_{i}))=h^{0}(C_{1},\mathcal{O}(a-a_{1})\oplus
\mathcal{O}(a-b_{1}))$.

\begin{enumerate} \item If $a<a_{1}$ and $a=b_{1}$,
condition \ref{condition} is still satisfied, but these points are
not in the closure $V_{d,a}$, (see \cite{Cos1}). \item If
$a=a_{1}$ and $a_{2}=0$, then there exists a family of surfaces
degenerating to a union of a cone of degree $b_{2}$ and a scroll
$X_{a_{1},b_{1}}$, so that these points are in the closure. \item
If $a=a_{1}$ and $a_{2}>0$, then by \cite{Cos1} the union
$X_{a_{1},b_{1}}\cup X_{a_{2},b_{2}}$ cannot be a flat limit of
scrolls.
\end{enumerate}

By degree we mean, the number of rational ruled surfaces
of this specifying type through $4d+1$ points, which we studied in
\cite{Mar1}.

For $a=\frac{d}{2}-1$, proposition \ref{dim} shows that
$R^{0}_{a}$ is a one codimensional subvariety. Its class
$(a-d+1)Y+D$ in the Chow ring $A(R^{0}_{d})$, is given by
Porteous formula in section \ref{porteous} and Lemma \ref{hypcl}.
Because the divisors in the Picard group $A^{1}(\overline{N})$
intersect the boundary transversally, it follows by
(\ref{generadores}) that this divisor extends in the moduli
compactification $\overline{M}$, to the divisor
$(a+d-1)\mathcal{L}_{i}+\mathcal{A}$. Therefore, the degree of the
stratum $\overline{M}_{\frac{d}{2}-1}$, that is, the degree of the
variety parametrizing rational ruled surfaces with a minimal
directrix of degree $\frac{d}{2}-1$, is the degree of the top
codimensional class $\left(
(a+d-1)\mathcal{L}_{i}+\mathcal{A}\right)\cdot
\mathcal{A}^{4d}\cdot \mathcal{L}_{1}\cdot \mathcal{L}_{2}\cdot
\mathcal{L}_{3}$ on $\overline{M}$, determined by general
$a-$planes $T_{a_{1}},\ldots,T_{a_{_{4d}}}$ and hyperplanes
$T_{3_{j}}$, $j=1,2,3$, (see \cite{Mar1}). Now since
$\mathcal{L}_{1}^{2}\cdot \mathcal{A}^{4d}=Q_{d}+Q_{d}^{b}$ in the
Chow ring $A(\overline{M}_{0,1}(\G,d))$, the result follows. \cqd

\subsubsection*{Remarks and Conclusions} Note that when we say
that $a$ takes the value $\frac{d}{2}-1$ it is implicitly
understood that $d$ is even, therefore we are only interested in
the numbers $Q^{b}_{d}$ listed in the table above for $d$ even.
The case $d=2$ corresponds to the cone
$\mathbb{P}(\mathcal{O}_{\mathbb{P}^{1}}\oplus
\mathcal{O}_{\mathbb{P}^{1}}(2))$. It is known that there are 4
cones through 8 general points. This is a degenerate case that was
excluded from the beginning (see 2.1 of \cite{Mar1}) for the
degree computation of the variety of rational ruled surfaces. As
an example we mention that according to Theorem \ref{degree},
there are 128054031870240 ruled surfaces with fixed parameters
$d=6$ and $a=2$. 

\end{document}